\documentclass[aps, nofootinbib,a4paper,floatfix]{revtex4-1}

\usepackage{graphicx,color}
\usepackage{amsmath,amssymb,amsfonts,amsthm,amscd,bm}
\usepackage{booktabs}
\usepackage{cancel}

\usepackage{amsmath}
\usepackage{amssymb}
\usepackage{amsfonts}
\usepackage{amsthm}
\usepackage{amscd}
\usepackage{bm}
\usepackage{graphicx}
\usepackage{booktabs}
\usepackage{listings}
\usepackage{xcolor,colortbl}
\usepackage{float}

\def\tilde{\widetilde}
\def\bar{\overline}
\def\hat{\widehat}
\def\*{\star}
\def\[{\left[}
\def\]{\right]}
\def\({\left(}      

\def\){\right)}

\def\frac#1#2{\dfrac{#1}{#2}}
\def\inv#1{\dfrac{1}{#1}}
\def\half{\tfrac{1}{2}}

\def\2pi{\hbox{$2\pi i$}}

\def\dsl{\raise.15ex\hbox{/}\kern-.57em\partial}
\def\Dsl{\,\raise.15ex\hbox{/}\mkern-.13.5mu D}
      \def\CF{{\cal F}}

   \def\CN{{\cal N}}

\def\2pi{\hbox{$2\pi i$}}

\def\dsl{\raise.15ex\hbox{/}\kern-.57em\partial}
\def\Dsl{\,\raise.15ex\hbox{/}\mkern-.13.5mu D}
\font\numbers=cmss12
\font\upright=cmu10 scaled\magstep1
\def\stroke{\vrule height8pt width0.4pt depth-0.1pt}
\def\topfleck{\vrule height8pt width0.5pt depth-5.9pt}
\def\botfleck{\vrule height2pt width0.5pt depth0.1pt}
\def\Zmath{\vcenter{\hbox{\numbers\rlap{\rlap{Z}\kern
    0.8pt\topfleck}\kern 2.2pt
    \rlap Z\kern 6pt\botfleck\kern 1pt}}}
\def\Qmath{
    \vcenter{\hbox{\upright\rlap{\rlap{Q}\kern3.8pt\stroke}\phantom{Q}}}}
\def\Nmath{\vcenter{\hbox{\upright\rlap{I}\kern 1.7pt N}}}
\def\Cmath{\vcenter{\hbox{\upright\rlap{\rlap{C}\kern
                   3.8pt\stroke}\phantom{C}}}}
\def\Rmath{\vcenter{\hbox{\upright\rlap{I}\kern 1.7pt R}}}
\def\Z{\ifmmode\Zmath\else$\Zmath$\fi}
\def\Q{\ifmmode\Qmath\else$\Qmath$\fi}
\def\N{\ifmmode\Nmath\else$\Nmath$\fi}
\def\C{\ifmmode\Cmath\else$\Cmath$\fi}
\def\R{\ifmmode\Rmath\else$\Rmath$\fi}
\def\barray{\begin{eqnarray}}
\def\earray{\end{eqnarray}}
\def\beq{\begin{equation}}
\def\eeq{\end{equation}}

\def\AA{\leavevmode\setbox0=\hbox{h}
\dimen0=\ht0 \advance\dimen0 by-1ex\rlap{\raise.67\dimen0\hbox{\char'27}}A}

\def\Arg{{\rm Arg}\,}

\def\kappa{c}

\def\Nc{N_{\rm max}} 

\makeatletter
\g@addto@macro\bfseries{\boldmath}
\makeatother

\makeatletter
\renewcommand*\env@matrix[1][\arraystretch]{%
  \edef\arraystretch{#1}%
  \hskip -\arraycolsep
  \let\@ifnextchar\new@ifnextchar
  \array{*\c@MaxMatrixCols c}}
\makeatother

\theoremstyle{plain}
\newtheorem{theorem}{Theorem}

\newtheorem{proposition}{Proposition}

\theoremstyle{remark}
\newtheorem{remark}{Remark}
\newtheorem{example}{Example}
\newtheorem{definition}{Definition} 
\newtheorem{proposal}{{\bf Proposal}}

\def\Q{\mathbb{Q}}
\def\Z{\mathbb{Z}}
\def\N{\mathbb{N}}
\def\nQ{\mathfrak{n}}
\def\rand{\mathfrak{r}}
\def\Shat{S}

\begin{document}

\title{
Riemann Hypothesis and  Random Walks:   the Zeta  case
}
\
\author{
 Andr\'e  LeClair\footnote{andre.leclair@gmail.com}
}
\affiliation{Cornell University, Physics Department, Ithaca, NY 14850} 

\begin{abstract}

In previous work it was shown 
 that if certain series  based on sums over primes of non-principal Dirichlet characters
have a conjectured random walk behavior,  
then the Euler product formula 
for its  $L$-function is valid to the right of the critical line $\Re (s) > \tfrac{1}{2}$, 
and the Riemann Hypothesis for this class of $L$-functions follows.  
Building on this  work,  here  we propose 
 how to extend this line of  reasoning  to  the Riemann zeta function and other principal Dirichlet 
 $L$-functions.   We apply these results to the study of the argument of the zeta function.  
In another   application,  we define and study  a 1-point correlation function of the Riemann zeros,    
    which  leads to the construction of a probabilistic  model for them. 
Based on these results we describe a new algorithm for computing very high Riemann zeros,  and 
we calculate the googol-th zero,  namely 
$10^{100}$-th zero to over 100  digits, far beyond what is currently known.    Of course use is made of the symmetry of the zeta function about the critical line.

\end{abstract}

\maketitle

\section{Introduction}

There are many generalizations  of Riemann's zeta function to other Dirichlet series,
which are also believed to satisfy a Riemann Hypothesis.       
A common opinion,   based largely on counterexamples,   is that the $L$-functions for which the  Riemann Hypothesis is true 
enjoy  both an Euler product formula and a functional equation.   However 
a direct connection between these properties and the Riemann Hypothesis has not been
formulated in a precise manner.       In \cite{EPF1,EPF2}  a concrete proposal making 
such a connection  was presented for Dirichlet $L$-functions,  and those based on cusp forms,
 due to  the validity of the Euler product formula to the right of the critical line.    
In contrast to the non-principal case, in this approach  the case of principal Dirichlet $L$-functions,  of which Riemann zeta is the simplest,  turned out to be  more delicate,  
 and consequently it was more  
difficult to state  precise results.      In the present work we  attempt to address further  this special case, although as we will explain the results are not as conclusive as for the non-principal case.     What is new that is presented here is a different way to understand the extent in which the truncated Euler product is a good approximation.
We then use this to approximate the argument of the zeta function on the critical line.   We also study 1-point statistics of the Riemann zeros,  in contrast to the 
2-point correlation functions that are widely studied.

Let $\chi (n)$ be a Dirichlet character modulo $k$ and $L (s, \chi)$ its $L$-function   with $s= \sigma + it$.
It satisfies the Euler product formula 
\beq
\label{Ldef} 
L (s, \chi) = \sum_{n=1}^\infty  \frac{\chi (n)}{n^s}  =  \prod_{n=1}^\infty  \(  1 -  \frac{\chi (p_n)}{p_n^s} \)^{-1} 
\eeq
where $p_n$ is the $n$-th prime.  
The above formula is   valid for $\Re (s) >1$  since  both sides converge absolutely. 
The  important distinction between principal verses non-principal characters is the following.  
  For non-principal characters the $L$-function has no pole at $s=1$, 
thus  there exists the possibility that the Euler product is valid  partway inside the strip,  i.e. has
abscissa of convergence $\sigma_c < 1$.   It was  proposed  in \cite{EPF1,EPF2} that 
$\sigma_c = \tfrac{1}{2}$ for this case.     In contrast,  now consider $L$-functions based on principal 
characters.   
The latter  character is defined as $\chi (n)=1$ if $n$ is coprime to $k$ and zero otherwise.
The Riemann zeta function is the trivial principal character of modulus $k=1$ with 
all $\chi (n) =1$.     $L$-functions based on principal characters  do have a pole at $s=1$, and  therefore
have abscissa of convergence $\sigma_c =1$, which  implies the Euler product  in the form  given above  strictly cannot be valid inside
the critical strip $0<\sigma < 1$.    Nevertheless,  in this paper we will show how  
a truncated  version of
the Euler product formula  can be approximately  valid for $\sigma > \tfrac{1}{2}$.     

The  primary aim of the work \cite{EPF1,EPF2}  was to determine  what specific properties of 
the prime numbers would imply that the Riemann Hypothesis is true.    This is the opposite 
of the more well-studied question of what the validity of the Riemann Hypothesis implies 
for the fluctuations in the distribution of primes.     The answer 
proposed  was simply based on the multiplicative independence of the primes,
which  to a  large extent  underlies their  pseudo-random behavior.    To be more specific,  
let $\chi (n) = e^{i \theta_n }$ for $\chi (n) \neq 0$.    
In \cite{EPF1,EPF2}  it was proven that if the series 
\beq
\label{BNt0} 
  B_N (t, \chi) =  \sum_{n=1}^N \cos \( t \log p_n + \theta_{p_n}  \) 
\eeq
is $O(\sqrt{N})$,    then the Euler product converges for $\sigma > \tfrac{1}{2}$  and the formula
\eqref{Ldef}  is valid to the right of the critical line.  
In fact,  we only need $B_N = O(\sqrt{N})$ up to logs (see Remark \ref{logs});  when we write 
 write $O(\sqrt{N})$,  it is implicit that this can be relaxed with logarithmic factors.    
   For non-principal characters the
allowed angles $\theta_n$ are equally spaced on the unit circle,  and it was conjectured  in
\cite{EPF2}  that the above series with $t=0$ behaves like a random walk  due to the 
multiplicative independence of the primes,  and this is the origin of the  
$O(\sqrt{N})$ growth.    Furthermore,  this result extends to all $t$ since domains of convergence
of Dirichlet series are always half-planes.    
    Taking the logarithm of \eqref{Ldef},  one sees that 
$\log L$ is never infinite to the right of the critical line and thus has no zeros there. 
This,  combined with the functional equation that relates $L(s)$ to $L(1-s)$,  implies there
are also no zeros to the left of the critical line,  so that all zeros are on the line.    
The same reasoning applies to cusp forms if one also uses a non-trivial result  of Deligne
\cite{EPF2}.

In this article we reconsider the principal Dirichlet case,  specializing to Riemann zeta itself 
since identical arguments apply to all other  principal cases with $k>1$ \footnote{In much more recent work a different approach to the non=principal case was studied based on 
M\"obius inversion \cite{ML}}.    
 Here all angles $\theta_n =0$, 
so one needs to consider the series
\beq
\label{BNt} 
  B_N (t) =  \sum_{n=1}^N \cos (t \log p_n ) 
\eeq
which now strongly depends on $t$.    
On the one hand,  whereas the case of principal Dirichlet $L$-functions is complicated by
the existence of the pole,  and,  as we will see,   one consequently  needs  to truncate the Euler product
to make sense of it,  on the
other hand  $B_N$ can be estimated using the prime number theorem since it does not involve
sums over non-trivial characters $\chi$,  and this aids the analysis.  
   This is in contrast to the non-principal case,  where,  however
well-motivated,   we had to
conjecture the random walk behavior alluded to above,  so in this respect  the principal case is  potentially simpler.    
To this end, 
a theorem of Kac (Theorem \ref{Kac} below)  nearly does the job:    $B_N (t) = O(\sqrt{N})$  
in the limit  $t \to \infty$,  which is also  a consequence of the multiplicative independence of the primes.
   This suggests  that one can also make sense of the Euler product formula in the limit 
 $t \to \infty$.
    However this is not enough for our  main purpose,  which  is to have
a similar result  for finite $t$ which we will develop.     

This article is mainly  based on our  previous work
 \cite{EPF1,EPF2}  but  provides a more detailed analysis and extends it in several ways.   
  It  was suggested in \cite{EPF1}  that one should truncate the 
series at an  $N$ that depends on $t$.     First, in the next section we explain how  a
simple  group structure
underlies a finite Euler product which relates it to a generalized Dirichlet series which is a 
subseries of the Riemann zeta function.   
Subsequently we  estimate  the error under truncation,  which shows explicitly   
 how this error is related to the pole at $s=1$,  as expected.    The remainder of the 
 paper, sections IV-VI,  presents various applications of these ideas.  We use them to study the argument of the zeta function.
   We present an algorithm to 
  calculate very high zeros,   far beyond what is currently known.       
 We also study the  statistical fluctuations of individual zeros,  in other words,  a 1-point correlation function. 
 
In many respects, our work is  related to   the work  of Gonek et. al. \cite{Gonek1,Gonek2},  which also considers a 
truncated Euler product.    The important difference is that  the starting point in \cite{Gonek1}  
is a hybrid version of the Euler product which involves both primes and zeros of zeta. 
Only after assuming the Riemann Hypothesis can one explain in that  approach  why the truncated product over
primes is a good approximation to zeta.
   In contrast,  here we do not assume anything about the zeros of zeta,   since the goal is to actually understand their location.           

We are unable to provide fully rigorous proofs of some of the statements below,  however we do provide supporting calculations and numerical work. 
In order to be clear on this,   below  ``{\bf Proposal}" signifies  the most important claims that we could not rigorously prove,  and should not be taken 
as a  ``{\bf Proposition}" in the usual formal mathematical sense.

\section{Algebraic  structure of finite Euler products}

\def\Q{\mathbb{Q}}
\def\Z{\mathbb{Z}}
\def\N{\mathbb{N}}
\def\nQ{\mathfrak{n}}

The aim of this section is  to define properly the objects we will be dealing with. 
In particular we will place   finite Euler products  on the  same footing  as 
other generalized Dirichlet series.   The results  are straightforward 
and are mainly definitions.  

\bigskip

\begin{definition}
\label{QNdef}  
Fix a positive integer $N$ and 
let $ \{p_1, p_2, \ldots p_N\}$  denote the first $N$ primes where $p_1 =2$.
   From this set one can  generate 
an abelian  group  $\Q_N$  of rank $N$ with  elements 
\beq
\label{QNeq}
\Q_N = \Bigl\{   p_1^{n_1} p_2^{n_2}   \cdots p_N^{n_N}, ~~n_i \in \Z ~\forall_i   \Bigr\}   
\eeq
where the group operation is ordinary multiplication.   
Clearly $\Q_N \subset \Q^+$ where $\Q^+$ are the positive  rational numbers.   
There are an infinite number of  integers in $\Q_N$  
 which form a  subset of the natural numbers $\N = \{ 1, 2, \ldots \}$.   We will denote this set as $\N_N \subset \N$, 
  and elements of this set simply as $\nQ$.
\end{definition}

\bigskip

\begin{definition}
Fix a positive integer $N$.   
For every integer $n \in \N$ we can define the character $c(n)$:
\barray
\nonumber
c(n)   &=& 1 ~~~~~{\rm if} ~ n \in \N_N \subset \Q_N
\\ 
\label{cn}
&=& 0 ~~~~~{\rm otherwise}   
\earray
Clearly,  for a prime $p$,   $c(p) = 0$  if $p> p_N$.
\end{definition}

\bigskip
\begin{definition}   Fix a positive integer $N$ and  let $s$ be a complex number.   
Based on $\Q_N$ we can define the infinite series 
\beq
\label{zetaN}
\zeta_N (s) =  \sum_{n=1}^\infty    \frac{c(n)}{n^s} =   \sum_{\nQ\ \in \N_N}   \inv{\nQ^s}  
\eeq
which is a generalized Dirichlet series.  
There are an infinite number of terms in the above series since $\N_N$ is infinite dimensional. 
\end{definition}

\bigskip
\begin{example}
For instance 
$$  \zeta_2 (s)  = 1 + \inv{2^s}  + \inv{3^s}  +  \inv{4^s} +   \inv{6^s}  +  \inv{8^s}  +  \inv{9^s}   + \inv{12^s} +  \ldots
$$
\end{example}

\bigskip

Because of the group structure of $\Q_N$,  $\zeta_N$  satisfies a finite  Euler product formula:
\begin{proposition}    
Let $\sigma_c$ be the abscissa of convergence of the series $\zeta_N (s)$  where 
$s = \sigma + i t$,  namely $\zeta_N (s)$ converges for $\Re (s) > \sigma_c$.    
Then in this region of convergence,  $\zeta_N$  satisfies a finite Euler product formula:
\beq
\label{zetaNEP}  
\zeta_N (s)   =   \prod_{n=1}^N   \( 1-  \inv{p_n^s} \)^{-1} 
\eeq
\begin{proof}   

\noindent
Based on the completely multiplicative property of the characters, 
\beq
\label{multc}
c(nm) = c (n) c(m)
\eeq
one has  
$$ \zeta_N (s) =   \prod_{n=1}^\infty   \( 1-  \frac{c(p_n)}{p_n^s} \)^{-1} $$
The result follows  then from the fact that   $c(p_n) = 0$  if  $n>N$.
\end{proof}
\end{proposition}     

\bigskip
\begin{example}
\label{geo}
Let $N=1$,  so that $\{ \nQ \} = \{1,2, 2^2, 2^3 \ldots \}$.     Then the above Euler product formula \eqref{zetaNEP}   is  simply the standard formula for the 
sum of a geometric series:
\beq
\label{geo2} 
\zeta_1 (s) =   \sum_{n=0}^\infty \inv{2^{ns}}  =   \inv{1 - 2^{-s} } 
\eeq
Here the abscissa of convergence is $\sigma_c = 0$.    
\end{example} 

\medskip

The series $\zeta_N (s)$  defined in \eqref{zetaN}  has some interesting properties:

\medskip

\noindent (i) ~  For finite $N$ the product is finite for $s \neq 0$,  thus the infinite series 
  $\zeta_N (s)$ converges 
for $\Re (s) > 0$ for any finite $N$.    

\medskip  

\noindent  (ii)    ~ Since the logarithm of the product is finite,   for finite $N$,   $\zeta_N (s)$ has no zeros nor poles 
for $\Re (s) > 0$.    Thus the Riemann zeros and the pole at $s=1$  arise from the primes
at infinity $p_\infty$, i.e.   
 in the limit $N \to \infty$.     
   In this limit all integers are included in the sum
 \eqref{zetaN} that defines $\zeta_N$ since $\N_\infty = \N$.   
 This is in accordance with  the fact that the pole 
 is a consequence of there being an infinite number of primes.   

\medskip

The property (ii) implies that,  in some sense,  the Riemann zeros condense out of the 
primes at infinity $p_\infty$.   
Formally one has 
\beq
\label{Ninf} 
\lim_{N \to \infty}  \zeta_N (s)   =   \zeta (s)
\eeq
However since $N$ is going to infinity,   the above is true only where the series  formally converges as a Dirichlet series,
which,  as discussed in the Introduction,  is $\Re (s) >1$.   Nevertheless,   for very large but finite $N$,   the function $\zeta_N$ can still  be a good approximation
to $\zeta (s)$ {\it inside the critical strip}  since for $N$ finite  there is convergence of
$\zeta_N (s)$ 
for $\Re (s) >0$.       
This is the subject of the next section,  where we show that 
a finite  Euler product formula 
is valid for $\Re (s) > \tfrac{1}{2}$ in a   manner that we will specify.

\section{Finite Euler product  formula at large $N$ to the right of the critical liine.}

In this section we propose that the Euler product formula  can be a very good approximation to $\zeta (s)$  for $\Re (s) > \tfrac{1}{2}$  and large $t$  if $N$ is chosen to depend on  $t$ 
in a specific way which was already proposed  in \cite{EPF1,EPF2}. 
  The new result presented here  is an  estimate  of the error due to the truncation.  

The random walk property we will build upon  is based on a central limit theorem of Kac \cite{Kac},  which largely 
follows from the multiplicative independence of the primes:

\begin{theorem} \label{Kac}
({\bf Kac})
Let  $u$ be a random variable uniformly distributed  on the interval $u\in [-T, T]$,   and define  the series 
\beq
\label{BN} 
  B_N (u) =  \sum_{n=1}^N \cos (u \log p_n ) .
\eeq
Then in the limit $N \to \infty$ and $ T \to \infty$,    $  B_N /\sqrt{N} $  approaches 
the normal distribution $\CN (0,1)$,  namely   
\beq
\label{prob}
\lim_{N \to \infty}  \lim_{T\to \infty}  P \Biggl\{  \frac{x_1}{\sqrt{2}}  <  \frac{ B_N (u)}{\sqrt{N} } < \frac{x_2}{\sqrt{2}}  \Biggr\}  = \inv{ \sqrt{2 \pi}} \int_{x_1} ^{x_2}
  e^{-x^2/2} dx
\eeq
where $P$ denotes the  probability for the set.    
\end{theorem}

We wish to use the above theorem to conclude something about $B_N (t)$ for 
a fixed,  non-random  $t$.    Based on Theorem \ref{Kac},   we suggest   the following
for  non-random,  but large  $t$:  
For any $\epsilon > 0$, 
\beq
\label{KacCor}
 \lim_{t \to  \infty}     B_N (t)  =    O(N^{1/2+\epsilon}).
\eeq
We could not rigorously  prove this statement,  however we can provide a heuristic argument.   As  $T \to \infty$,   even though $u$ is random,   the vast majority of them  
are tending to $\infty$.  One then uses the normal distribution in Theorem \ref{Kac}.   
In the following we will provide indirect numerical evidence.

 \begin{remark}
 \label{logs} 
  The proof of convergence of the Euler product  in \cite{EPF2}  is not spoiled if the bound on $B_N$ is relaxed up to logs.  For instance, if in the limit
  $t \to \infty$, 
 $B_N = O(\sqrt{N \log \log N})$,  as suggested by the law of iterated logarithms relevant to central limit theorems,  this is fine,  as is $B_N = O( \sqrt{N} \log^a N)$ for any positive power 
 $a$.    
 \end{remark}

 A consequence of  Theorem \ref{Kac}  and the comments following it  is that the Euler product formula is valid  to the right of the
 critical line in the limit
  $t \to \infty$,  at least formally.  Namely
 for $\sigma > \tfrac{1}{2}$,   
\beq
\label{EPF}
\lim_{t\to \infty}   \zeta (\sigma + i t)  = \lim_{N \to \infty}  \lim_{t\to \infty}
\prod_{n=1}^N  \( 1- \inv{p_n^{\sigma + i t }} \)^{-1}  
\eeq
As shown  in \cite{EPF1,EPF2} and discussed in the Introduction,  this  formally follows from the $\sqrt{N}$ growth of $B_N$.   
The problem with the above formula is that due to the double limit on the RHS,  it is not rigorously defined.
For instance,  it could depend on the order of limits.   
It is thus desirable to have a version of \eqref{EPF}
 where $N$ and $t$ are taken to infinity simultaneously.
 Namely,  we wish to truncate the product  at  an   $N (t)$ that depends
 on $t$ with the property that $\lim_{t\to \infty} N (t) = \infty$.     One can then replace
 the double limit on the RHS of \eqref{EPF}  with one limit $t \to \infty$,   or equivalently
 $N (t) \to \infty$.
 
   There is no unique choice for $N (t)$,  but there is 
 an optimal  upper limit,  $N(t) <  \Nc (t)  \equiv [t^2]$,  with $[t^2]$ its integer part,   which we now describe.  
  We can use the prime number theorem to estimate $B_N (t)$:
 \barray
 \label{BNpnt} 
 B_N (t)  \approx \int_2^{p_N}  \frac{dx}{\log x}   \cos ( t\log x )   
 &=& \Re\(   {\rm Ei} \( (1 + it) \log p_N \)   \)  
 \\ \nonumber
& \approx&  \frac{p_N}{\log p_N}  \( \frac{t}{1+t^2 } \)  \sin \( t \log p_N \)  
 \earray
 where  ${\rm Ei}$ is the usual  exponential-integral function,  and we have used
 \beq
 \label{Eiasym}
 {\rm Ei} (z) = \frac{e^z}{z} \( 1 + O \( \frac{1}{z}\) \)
 \eeq
 The prime number theorem implies $p_N \approx N \log N$.   Using this in
 \eqref{BNpnt} and  imposing  $B_N (t) < \sqrt{N}$ leads to $N < [t^2]$.  
 
 Based on the above,  henceforth we will always assume the following properties of $N(t)$:
 \beq
 \label{Ntprop}
 N(t) \leq   \Nc (t) \equiv [t^2]~~~~{\rm with} ~~~~~\lim_{t\to \infty}   N(t) = \infty,
 \eeq
 and will not always display the $t$ dependence of $N$.    
Equation \eqref{EPF} now  formally becomes
 \beq
 \label{EPFa}
 \lim_{t \to \infty}   \zeta (s) = \lim_{t \to \infty}  \prod_{n=1}^{N (t)}  \( 1- \inv{p_n^s} \)^{-1} , ~~~~~
 {\rm for} ~ \Re (s) > \tfrac{1}{2}
 \eeq
 \noindent
 Extensive and compelling numerical evidence supporting the above formula was
 already presented in \cite{EPF1}.

Based on the above results we are now in a position to study the following 
important question.
If we fix  a finite  but large 
 $t$,  and truncate the Euler product at $N(t)$,  which is finite,  what is the error in
the approximation to $\zeta$ to the right of the critical line?     We estimate 
this error as follows:   
 \begin{proposal}\label{EPFcut}
 Let $N=N(t)$ satisfy \eqref{Ntprop}.   Then 
 for $\Re (s) > \tfrac{1}{2}$ and large $t$, 
 \beq
 \label{EPFNc}
 \zeta (s) = \prod_{n=1}^{N(t) }  \(  1 - \inv{p_n^s} \)^{-1}   \exp \( R_{N} (s) \) 
 \eeq
 where  $\zeta (s)$ is the actual $\zeta$ function defined by analytic continuation and 
 \beq
 \label{rNcapprox}
 R_{N} (s)  =    \inv{(s-1)}  \,  O\(  \frac{ N^{1-s} }{\log^s N }  \).
 \eeq
 $R_N$  
  is finite (except at the pole $s=1$)  and 
 satisfies 
 \beq
 \label{RNlim}
 \lim_{t\to \infty}  R_{N(t)} (s) = 0,
 \eeq 
 namely the error goes to zero as $t \to \infty$.  
 \end{proposal}
 
We provide the following supporting argument,  although not a rigorous proof,   for this Proposal.      
From  \eqref{EPFa},    one concludes that \eqref{EPFNc} must hold in the limit of large $t$ 
with  $R_N$ satisfying \eqref{RNlim}.   The logarithm of \eqref{EPFNc}
reads
\beq
\label{logEPF}
\log \zeta (s) =  - \sum_{n=1}^{N }  \log \( 1 - \inv{p_n^s} \)  +   R_{N} (s) 
\eeq
First assume $\Re (s) > 1$.  Then 
in the limit of large $t$,   the error upon truncation is the part that is neglected 
in \eqref{EPFa}:
\beq
R_N (s)  = -  \sum_{n=N+1}^\infty \log \( 1 - \inv{p_n^s} \)
\eeq
Expanding out the logarithm,  one has  
\barray
\nonumber 
R_{N} (s)  &\approx & \sum_{n=N }^\infty  \inv{p_n^s}    
\\ 
\label{rNO1}
&\approx&  \int_{p_{N}}^\infty \frac{dx}{\log x}   \inv{x^s}  
\approx   \inv{(s-1)} \frac{p_{N}^{1-s}}{\log p_{N}}  
\earray
where in the second line we again used the prime number theorem to approximate the sum over primes.  
Next using $p_N \approx N \log N$,   one obtains \eqref{rNcapprox}. 
Finally,  the above expression can be continued into the strip $\sigma > \tfrac{1}{2}$ if
$N(t) < [t^2]$  since 
$N(t)^{1-s}/t  <  N^{1/2 - s}$ which goes to zero as 
$N \to \infty$ 
  if $\Re(s) > \tfrac{1}{2}$.     The latter also implies \eqref{RNlim}.

\begin{figure}[t]
\centering\includegraphics[width=.5\textwidth]{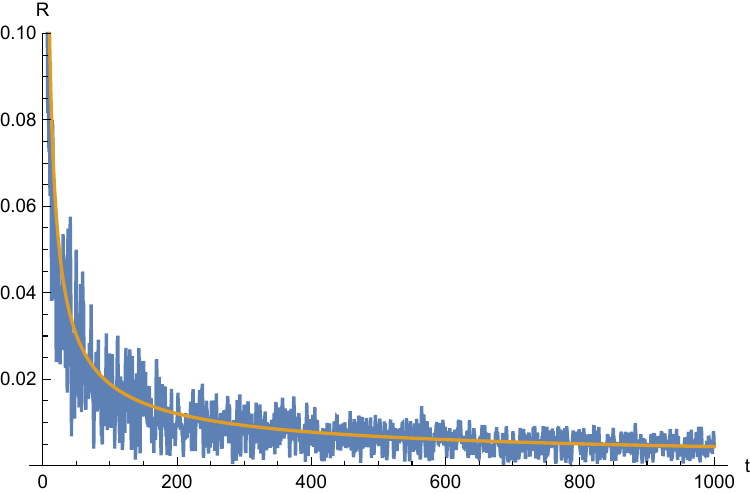}
\caption{The error term $|R_{N} (s)|$  with $N(t) = \Nc (t) = [t^2]$ for $\Re (s) = 3/4$ inside the critical strip as a function of
$t$.   The fluctuating (blue) curve is $|R_{N}|$  computed directly from the definition 
\eqref{EPFNc} with $\zeta (s)$ the usual analytic continuation into the strip.      The smooth (yellow) curve  is the approximation 
$ R_{N} (s)  =    \inv{(s-1)}  \,   \frac{ N^{1-s} }{\log^s N }  $ based on 
\eqref{rNcapprox}.   
  }
\label{RNc}
\end{figure}

 Proposal  \ref{EPFcut}  makes it clear  that the need for a cut-off $N< \Nc$ originates from 
 the pole at $s=1$,  since as long as $s\neq 1$,  the error $R_{N} (s)$ in \eqref{rNcapprox}
  is finite.    The error becomes smaller and smaller the further one is from the pole,  i.e. 
  as $t \to \infty$.    
In Figure \ref{RNc} we  numerically illustrate Proposal \ref{EPFcut} inside the critical strip.  

\begin{remark}
For estimating errors at large $t$ the following formula is   useful:  
\beq
\label{RNO}
|R_{N(t)} (s) | \sim  \frac{  N(t) ^{1-\sigma}  }{t} \sim   \inv{t^{2\sigma -1}}
\eeq
\end{remark}

\begin{proposal}
\label{RH} 
Assuming  Proposal  \ref{EPFcut},   all non-trivial  zeros of $\zeta (s)$  
are on the critical line.  
\begin{proof}
Taking the logarithm of the truncated Euler product,  one obtains 
\eqref{logEPF}.    If there were a zero  $\rho$ with $\Re (\rho ) > \tfrac{1}{2}$,  then 
$\log \zeta (\rho ) = - \infty$.    However the right hand side of \eqref{logEPF}  is 
always finite,  thus there are no zeros to the right of the critical line.   The functional
equation relating $\zeta (s) $ to $\zeta (1-s)$ shows there are also no zeros to the left
of the critical line.    
\end{proof}
\end{proposal}
 
  \begin{remark} 
 Interestingly,   Proposal  \ref{EPFcut}  and Theorem  \ref{RH} 
 imply  that proving the  validity of the Riemann Hypothesis is under better 
 control the higher one moves up the critical line.  For instance,  it is known that
 all zeros are on the line up to $t\sim 10^{13}$,  and beyond this,   the error $R_N$ is
 too small to spoil the validity of the Riemann Hypothesis.    Henceforth,  we assume the RH. 
 \end{remark}

 \section{$1$-point correlation function of  the Riemann zeros}

\label{1point_sec} 
 
\def\tnN{{t_{n;N}}}

\def\Fluc{\CF}  

Montgomery conjectured that the pair correlation function of ordinates of the Riemann zeros
on the critical line satisfy GUE statistics \cite{Montgomery}.       Being a 2-point correlation function,
it is a reasonably 
complicated statistic.   In this section we propose a simpler  1-point correlation function 
that captures the statistical fluctuations of individual zeros.

Let $t_n$ be the exact  ordinate of the $n$-th  zero on the critical line,  with $t_1 = 14.1347...$ and so 
forth.     The single equation $\zeta (\rho) =0$ is known to have  an infinite number of  non-trivial solutions 
$\rho = \tfrac{1}{2} + i t_n$.     In \cite{Trans},  by placing the zeros in one-to-one correspondence
with the zeros of a cosine function,   the single equation  $\zeta (\rho) =0$ was replaced by
an infinite number of equations,  one for each $t_n$ that depends only on $n$: 
\beq
\label{transeq}
\vartheta (t_n) + \lim_{\delta \to 0^+}   \arg \zeta (\tfrac{1}{2} + \delta + i t_n ) = (n-\tfrac{3}{2} ) \pi 
\eeq
where $\vartheta$ is the Riemann-Siegel function:
\beq
\label{RS}
\vartheta (t) =  \Im \log \Gamma( \tfrac{1}{4} +  \tfrac{i t}{2} )  -  t \log \sqrt{\pi}.
\eeq
The equation \eqref{transeq} involves the important function 
\beq
\label{SofT}
S(t) = \inv{\pi}  \lim_{\delta \to 0^+} \arg \, \zeta (\half + \delta + it  ).
\eeq
It is important that the $\delta \to 0^+$  approaches the critical line from the right,  since this is where
the Euler product formula is valid in the sense described above.  
This equation was used to calculate zeros very accurately in \cite{Trans},  up to thousands of
digits.    There is no need for a cut-off $\Nc$ in the above equation since 
the $\arg \zeta$ term is defined for arbitrarily high $t$ by 
 standard   analytic continuation.    
One aspect of this equation is the  following theorem: 

\begin{theorem} 
\label{FLthm}   ({\bf Fran\c ca-LeClair})   If there is a unique solution to the 
equation \eqref{transeq}  for every positive integer  $n$,  then the Riemann Hypothesis is true,  and 
furthermore,  all zeros are simple.  
\end{theorem}  

\begin{remark} 
Details of the proof are in \cite{Trans}.    
   The main idea is that if there is a unique solution,  then the zeros 
are enumerated by the integer $n$ and can be counted along the critical line,  and the resulting counting formula  coincides   with a well known result 
  due to Backlund for the number of zeros 
in the {\it entire}  critical strip.    The zeros are simple because the zeros of the cosine are simple.     
The above theorem is another approach towards proving the Riemann Hypothesis,    however it is not entirely independent of 
the above approach based on the Euler product formula,  in particular Theorem \ref{RH}.   
In \cite{Trans},  we were unable to prove there is a unique solution because we did not  have sufficient control over   the 
relevant properties of the  argument of $\zeta$ on the critical line. 
   \end{remark}

If the $\arg \zeta$ term is ignored,  then there is indeed a unique solution 
for all $n$ since $\vartheta (t)$  is a monotonically increasing function of $t$.   
Using its asymptotic expansion for large $t$,  equation  \eqref{RSseries} below, and dropping the $O(1/t)$ term,   then the 
solution is 
\beq
\label{tLambert} 
\tilde{t}_n =  \frac{2 \pi ( n - \tfrac{11}{8} )}{W \( (n-\tfrac{11}{8})/e \) }
\eeq
where $W$ is the Lambert $W$-function.   
The only way there would fail  to be a solution is if $\Shat (t)$ is not well defined 
for all $t$.      We point out that the Lambert function was used in connection with the Riemann zeros in 
\cite{Riguidel},  however the meaning does not seem to be the same as in this article.

\begin{figure}[t]
\centering\includegraphics[width=.5\textwidth]{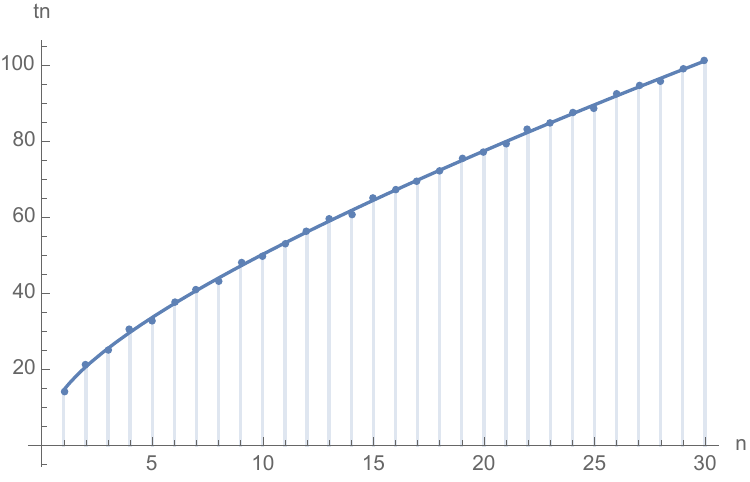}
\caption{ The first 30 Riemann zeros $t_n$.   The smooth curve is the approximation 
$\tilde{t}_n$ in \eqref{tLambert},  whereas the dots are the actual zeros $t_n$. }
\label{tn}
\end{figure}

The fluctuations in the zeros  come from $\arg \zeta$ since $\tilde{t}_n$ is a smooth
function of $n$.   These small fluctuations  are shown in  Figure \ref{tn}.   Let us  define
  $\delta t_n =  t_n -  \tilde{t}_n$.
One needs to properly normalize $\delta t_n$, 
taking into account that the spacing between zeros decreases as $2\pi /\log n$.   
To this end we expand the equation \eqref{transeq}  around $\tilde{t}_n$.   
Using $\vartheta (\tilde{t}_n ) \approx (n- \tfrac{3}{2})\pi$,  one obtains 
$\delta t_n \approx - \pi \Shat (t_n)/ \vartheta' (\tilde{t}_n) $ where $\vartheta' (t)$ is the derivative
with respect to $t$.    Using  $\vartheta' (t) \approx \tfrac{1}{2} \log (t/2 \pi e)$,   this 
 leads us to define
\beq
\label{deltadef} 
\delta_n \equiv  \frac{ (t_n - \tilde{t}_n )} {2 \pi}   \log \(  \frac{\tilde{t}_n}{2 \pi e} \)   \approx -  \Shat (t_n)
\eeq
The  probability distribution of the set 
\beq 
\label{DeltaM}
\Delta_M  \equiv   \Bigl\{ \delta_1, \delta_2,  \ldots , \delta_M \Bigr\}
\eeq
  for large $M$ is then an interesting property to study.  Here ``probability" is defined as frequency of occurrence.   
 The origin of the statistical fluctuations of $\Delta_M$ is  ultimately the fluctuations in the primes.

     \def\sigmag{\sigma_1} 
  
\begin{figure}[t]
\centering\includegraphics[width=.5\textwidth]{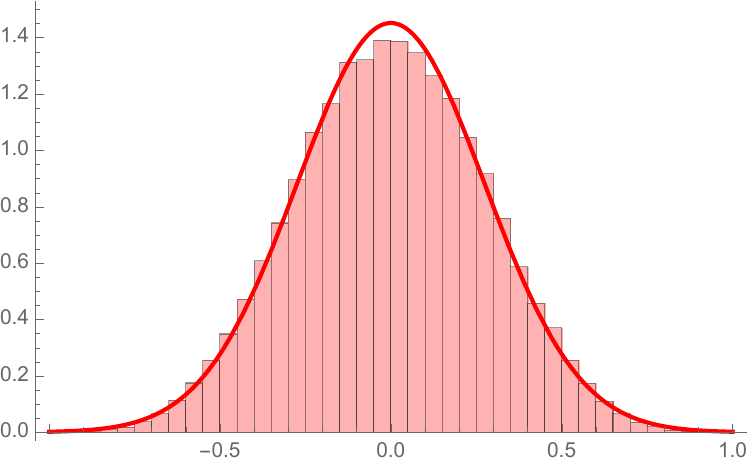}
\caption{The probability distribution for the set $\Delta_M$  defined in \eqref{DeltaM} 
for $M=10^5$.   
The smooth curve is the normal distribution $\CN(0,\sigmag)$ with 
$\sigmag= 0.274$.  }
\label{NormalZeros}
\end{figure}

  In Figure \ref{NormalZeros} we plot the distribution of $\Delta_M$ for $M=10^5$.  
  It closely resembles a normal distribution.
  Let us  suppose  $\Delta_M$   does indeed satisfy  a normal distribution 
 $\CN(\mu, \sigmag)$.   
 Using  some known properties of $\Shat (t_n)$,   together with 
 the equation \eqref{deltadef},   we can propose then 
 the following.   
First,  one expects that the average of $\delta_n$ is zero since it is known that the
average of $\Shat (t)$ is zero,  thus $\mu=0$.        Up to the height $t$ that we have studied, 
 $\Shat (t)$ is nearly always on the principal branch,  i.e. $-1 < S(t) < 1$ up to some reasonably high $t$ on the order of $t=10^6$ or more.  
Then at each jump by $1$ at $t_n$,   on average $\Shat (t_n)$ passes through zero.
This implies that the average $\bar{ | \Shat (t_n) |} \approx 1/4$.   
For a normal distribution $\bar{ | \Shat (t_n) |} =\sqrt{\tfrac{2}{\pi}}  \, \sigmag$.   
Thus one expects  the standard deviation  $\sigmag$ of $\Delta_M$ to be 
$\sigmag \approx \sqrt{\pi/32} = 0.313..$.   
   In Figure \ref{NormalZeros} 
we present results for the first $10^5$-th known  exact zeros.    The distribution function 
 fits  a normal distribution  with $\sigmag = \sqrt{\pi/32}$  rather  well.   
 Performing a fit,  one finds $\sigmag \approx 0.27$.   For higher values of $M$ around $10^6$,  a fit gives
 $\sigmag \approx 0.3$,  which is closer to the predicted value.     We emphasize however that this  approximate prediction for $\sigma_1$ assumes $S(t)$ is on the principal branch, 
 which  is not expected to hold for arbitrarily high $t$.

\bigskip

\def\that{\hat{t}}

If  we approximate the distribution of $\Delta_M$ as normal,  then  we  can construct a simple  probabilistic model of the
Riemann zeros:     

\bigskip

\begin{definition}  {\bf ~ A probabilistic model of the Riemann zeros.}  ~
\label{randommodel}
 Let $\rand$ be a random variable with normal distribution
$\CN (0,\sigmag)$. 
Then a probabilistic  model of the  zeros $t_n$ can be defined as 
the set $\{ \that_n \}$,  where 
\beq
\label{tRand}
\that_n   \equiv   \tilde{t}_n +   \frac{2 \pi\, \rand }{\log ( \tilde{t}_n / 2\pi e )}  
\eeq
and $\tilde{t}_n$ is defined in \eqref{tLambert}.   In the above formula $\rand$ is chosen at random independently for each $n$.  
\end{definition} 

\bigskip

The statistical model \eqref{tRand} is rather simplistic since 
it is just  based on a normal distribution for $\rand$ and $\tilde{t}_n$  is smooth and
completely deterministic.     A natural  question then 
arises.    Does the pair correlation function of $\{ \that_n \}$ satisfy GUE statistics as does 
the actual zeros $\{ t_n \}$?       
It  is certainly  interesting to 
study the 2-point correlation function of    $\{ \that_n \}$.   
Montgomery's pair correlation conjecture can be stated as follows.
Let $\CN(T)$ denote the number of zeros up to height $T$,  where
$\CN(T) \approx  \tfrac{T}{2\pi}  \log  \( \tfrac{T}{2 \pi e}   \)$.    Let $t,t'$ denote zeros in the range
$[0,T]$.     Then in the limit of large $T$:  
\beq
\label{mont}
\inv{\CN (T)}  \sum_{\alpha < d(t,t') < \beta }   1   \sim  \int_\alpha^\beta  du 
\(  1 - \frac{\sin^2 (\pi u)}{\pi^2 u^2 }  \)  
\eeq
where $d(t,t')$ is a normalized distance between zeros
$d(t,t')  =  \tfrac{1}{2\pi} \log \( \tfrac{T}{2 \pi e} \) (t-t')  $.    

\begin{figure}[t]
\centering\includegraphics[width=.5\textwidth]{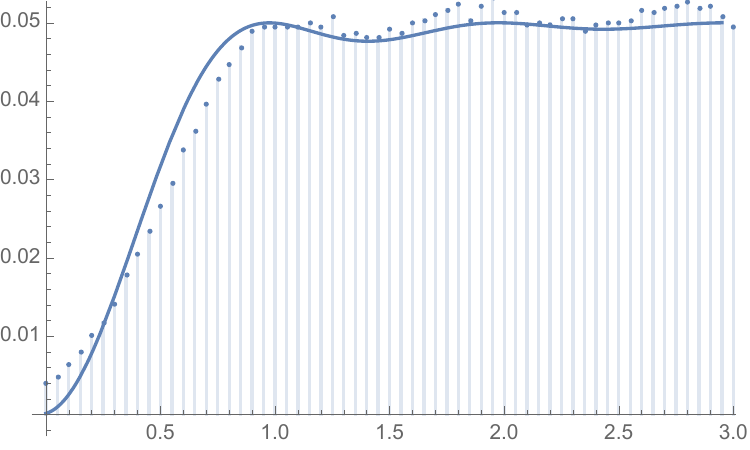}
\caption{The pair correction function of $\{ \that_n \}$  defined in \eqref{tRand} 
for $n$ up to $10^5$ where 
the standard deviation of $\rand$  was taken to be  $\sigmag = 0.274$.   The solid curve is the 
GUE prediction.   The parameters in \eqref{mont}  are $\beta = \alpha + 0.05$ with $\alpha = 
(0, 0.05, 0.10, \ldots, 3)$ and the $x$-axis is given by $x=(\alpha + \beta)/2$.   }
\label{pair}
\end{figure}

In Figure \ref{pair}   we   plot the pair correlation  function for the first $10^5$-th 
$\that_n$'s.    We chose  $\sigmag = 0.274$  since in this range of $n$ 
this gives a better fit to the normal distribution of the 1-point function.  
 The results  are  reasonably  close to the GUE prediction \eqref{mont},  especially considering
 that for just the first $10^5$ true zeros the fit to the GUE prediction is not perfect;    for much higher zeros
 it is significantly  better \cite{OdlyzkoPair}.

 \section{Computing very high zeros  from the primes}

 This section can be viewed as providing  additional  numerical evidence for some  of the previous results.   
We will be calculating $\Shat (t)$ from the primes using the truncated Euler product.
 Since this  requires $\Re (s) \to \tfrac{1}{2}^+$,  this is pushing
the limit of the validity of the truncated  Euler product formula,   nevertheless we will obtain reasonable results. 
We emphasize that this method has nothing to do with the   random model for the zeros in Definition \ref{randommodel},
but rather relies on the Euler product formula to calculate $\Shat (t)$.    

Many very high zeros of $\zeta$ have been computed numerically,   beginning with the work
of Odlyzko.     All zeros up to the $10^{13}$-th  have been computed and are all on the
critical line \cite{Gourdon}.     Beyond this the computation of zeros 
remains a challenging open problem.
However some zeros around the $10^{21}$-st and $10^{22}$-nd  are known \cite{Odlyzko}.     
In this section we describe a new and simple algorithm for computing very high zeros 
based on the above reasoning.   It will allow us to go much higher than the known 
zeros since it does not require  numerical implementation of the $\zeta$ function itself,  but rather 
only requires knowledge of some of the lower primes.

Let us first discuss the numerical challenges involved in computing  high zeros  from the 
equation \eqref{transeq}  based on the 
standard Mathematica package.    
The main difficulty is that one needs to implement the $\arg \zeta$ term.   
Mathematica computes  $\Arg \zeta$,  i.e. on the principal branch,   however near a zero
this is likely to be valid based on the discussion in section IV.    The main problem is that  Mathematica 
can only compute $\zeta$ for $t$ below some  maximum value around 
$t = 10^{10}$.    This was sufficient to calculate up to the $n=10^9$-th zero 
from \eqref{transeq} in \cite{Trans}.       The $\log \Gamma$  term must also be implemented to very high $t$,
which is also limited in  Mathematica.

We deal with these difficulties first by computing $\arg \zeta$ from the Euler product  formula 
 involving a finite sum over primes.          Then,  the $\log \Gamma$ term can be accurately computed 
using corrections to Stirling's formula:
\beq
\label{RSseries} 
\vartheta (t) =  \frac{t}{2} \log \( \frac{t}{2\pi e} \) - \frac{\pi}{8}  + \inv{48\,t}   
 + O(1/t^3) 
\eeq

Let $\tnN$ denote the ordinate of the $n$-th zero computed using the first $N$ primes based on \eqref{transeq}.  
For high zeros,   it is approximately the solution to the following equation
\beq
\label{EqtoSolve}
 \frac{\tnN}{2} \log \( \frac{\tnN}{2\pi e} \) - \frac{\pi}{8}   
    - \lim_{\delta \to 0^+} \Im \sum_{k=1}^N \log \( 1- \inv{p_k^{1/2 + \delta + i \tnN}} \) 
=  (n-\tfrac{3}{2} ) \pi 
\eeq
where it is implicit that $N < \Nc (t) =[t^2]$. 
The important property of this equation is that it no longer makes any reference to 
$\zeta$ itself.    It is straightforward to solve the above equation with standard 
root-finder software,  such as FindRoot in Mathematica.      

One can view the computation of $t_n$ as a kind of Markov process.       
If one includes no primes,  i.e. $N=0$,   and drops the next to leading $1/t$ corrections, 
then the solution is unique and explicitly given by $t_{n;0} = \tilde{t}_n$  in  terms of the Lambert $W$-function in 
\eqref{tLambert}.  
    One  then goes from $t_{n;0}$ to $t_{n;1}$ 
by finding the root to the equation for $t_{n;1}$ in the vicinity of $t_{n;0}$,  then similarly 
$t_{n;2}$ is calculated based on $t_{n;1}$ and so forth.   
At each step in the process one includes one additional prime,  and this slowly approaches 
$t_n$,  so long as $N(t) < \Nc (t)$.    In practice we did not follow this iterative procedure,  but rather fixed $N$ and simply
solved \eqref{EqtoSolve} in the vicinity of $\tilde{t}_n$.  

We can  estimate the error in computing the zero $t_n$ from the primes using equation \eqref{EqtoSolve}
as follows.    As in Section \ref{1point_sec},   we expand the equation \eqref{transeq} now around $t_{n;N}$
rather than $\tilde{t}_n$.    
One obtains 
$$t_n - \tnN = - \pi \, d\Shat _N /  \vartheta' (\tnN)$$
 where $d\Shat_N$ is the error in computing $\Shat (t)$ 
from the primes.       Using \eqref{rNO1},  we have 
$$d\Shat_N   =  \tfrac{1}{\pi}   \Im R_N (s=\tfrac{1}{2} +i t )  \approx \frac{ \sqrt{p_N}}{\pi t \log p_N}   \cos (t \log p_N).$$    
Now from the prime number theorem,  $p_N \approx N \log N$.     Recall $N$ is cut off at $\Nc = [t^2]$,  which cancels
the $1/t$ in the previous formula.     Finally it is meaningful to normalize the error by the mean spacing $2\pi/ \log n$.   
The result is 
\beq
\label{error}  
\frac{  t_n -  \tnN }{\scriptstyle{2\pi/ \log n}}    \approx  \inv{\pi \sqrt{\log N}}  \cos \(  t_n \log p_N \)  
\eeq
where we have used $\tnN \approx  \tilde{t}_n \approx 2 \pi n/\log n$.   
The left  hand side represents the ratio of the error to the mean spacing between zeros at that height.   
Again, it is implicit that $N < [t_n^2]$.     The interesting aspect of the above formula is 
that the relative error decreases with $N$,   although rather slowly.    The cosine factor also implies there are 
large scale oscillations around the actual $t_n$.        

For very high $t$,   $\Nc (t) =[t^2]$ is extremely large and it is not possible in practice 
to work with  such a large number of primes.    This  is the primary  limitation 
to the accuracy we can obtain.     We will limit ourselves to the relatively small  $N=5 \times 10^6$ primes.
Let us verify the method by comparing with some  known zeros around $n=10^{21}$ and $10^{22}$. 
     The results are shown in Table \ref{1021}.     Equation \eqref{error} predicts  
     $t_n -  t_{n;N} \approx 0.01$ for these  $n$ and $N$,   and inspection of the table shows this is a good estimate.    
         Odlyzko was of course able to calculate
 more digits;  our  accuracy 
can be improved by increasing $N$ in principle.           We also checked some zeros around the $n=10^{33}$-rd 
computed by Hiary \cite{Hiary},  again with favorable results.     

\begin{table}
\begin{center}
\begin{tabular}{|c|c|c|}
\hline\hline
$n$&  ~ ~~$\tnN$~~ ~ &  $t_n$ (Odlyzko) \\
\hline\hline 
$10^{21} -1 $  & $144176897509546973538.205$  &$\sim   .225 $\\
$10^{21} $  & $144176897509546973538.301$  &$\sim   .291 $\\
$10^{21} +1$  & $144176897509546973538.505$   & $\sim .498$\\
$\vert$ & $\vert$ & $\vert$ \\
$10^{22} -1$  & $1370919909931995308226.498$ & $\sim.490$ \\
$10^{22}  $  & $1370919909931995308226.614$ & $\sim.627$ \\
$10^{22} +1$  & $1370919909931995308226.692$ &$ \sim.680$ \\
\hline\hline 
\end{tabular}
\end{center}
\caption{Zeros around the $n=10^{21}$-st and $10^{22}$-nd   computed from \eqref{EqtoSolve} 
with $N=5\times 10^6$ primes.  We fixed $\delta = 10^{-6}$.  Above,
 $\sim$ denotes the integer part of
the second column.}
\label{1021}
\end{table}

Having made this check,  
let us now go far beyond this and compute the $n=10^{100}$-th zero  by the same method. 
Again   using  only  $N=5 \times 10^{6}$ primes,  we found the
following $t_n$:  
\barray
\nonumber 
n= 10^{100} {\rm -th ~ zero}: &&
\\ 
\nonumber 
  && {t_n= 280690383842894069903195445838256400084548030162846}\atop 
 {~~~~~~~~~~~~~~045192360059224930922349073043060335653109252473.244....}  
\earray
Obtaining this number took only a few minutes on a laptop using Mathematica.   
 We are confident that the last $3$ digits $\sim .244$ are  correct  since   we checked that 
 they didn't 
 change between   $N=10^6$ and $5\times 10^6$.    Furthermore,   $3$ digits is consistent with  \eqref{error},   which predicts 
 that for these $n$ and $N$,   $t_n - \tnN \approx 0.002$.    
      We calculated the next zero to be $\sim .273$.    
      
  We were able to extend this calculation to the $10^{1000}$-th zero without much difficulty.     As equation \eqref{error} shows,   the relative error 
  only decreases as one increases $t$.    
  It is also  straightforward  to extend this method to all primitive Dirichlet $L$-functions and those based on cusp forms  using the transcendental equations in 
  \cite{Trans}  and  the results in \cite{EPF2}.

\section*{Acknowledgments}

We wish to thank Denis Bernard,   Guilherme Fran\c ca,  Ghaith Hiary,  Giuseppe Mussardo,   and German Sierra  for discussions.   We also thank  the Isaac Newton Institute for 
Mathematical Sciences for their hospitality in the final stages of this work (January 2016).


\begin{thebibliography}{99}


\bibitem{EPF1}    G.  Fran\c ca and A. LeClair, 
``On the validity of the Euler product inside the critical strip", 
arXiv:1410.3520 [math.NT]. 

\bibitem{EPF2}    G.  Fran\c ca and A. LeClair, 
``Some Riemann Hypotheses from Random Walks over Primes", 
Communications in Contemporary Mathematics (2017) 1750085, 
arXiv:1509.03643  [math.NT]. 

\bibitem{ML}
G. Mussardo and A. LeClair,
{\it 
Randomness of M\" obius coefficents and brownian motion: growth of the Mertens function and the Riemann Hypothesis}
arXiv:2101.10336 [math-NT],  to appear in JSTAT.  



\bibitem{Kac}    M. Kac, 
{\it Statistical Independence in Probability,  Analysis and Number Theory},
The Mathematical Association of America,  New Jersey, 1959. 

\bibitem{Gonek1}   S.  M.  Gonek,  C.  P.  Hughes,  and J. P.  Keating,  
Duke Math. J {\bf 136} (2007) 507.

\bibitem{Gonek2}   S. M. Gonek,  Trans.  Amer.  Math.  Soc.   {\bf 364}  (2011)  2157.  


\bibitem{Trans}    G.  Fran\c ca and A. LeClair,
{\it Transcendental equations satisfied by the individual zeros of Riemann $\zeta$,  Dirichlet, and modular 
$L$-functions}, 
Commun.  Number Theory and Phys.   {\bf 9} (2015) 1-50.   


\bibitem{Riguidel}  M.  Riguidel, 
{\it The Two-Layer Hierarchical Distribution Model of Zeros of Riemann?s Zeta Function along the Critical Line},  Information 2021, 12, 22. https://doi.org/10.3390/info12010022.


\bibitem{Montgomery}   H. Montgomery, in 
{\it Analytic number theory,  Proc.  Sympos. Pure Math.  XXIV }   
(Providence,  RI:  AMS, 1973).   

\bibitem{OdlyzkoPair}   A. M. Odlyzko,  Math.  Comp. {\bf 48}.  273  (1987).  


\bibitem{Gourdon} 
X. Gourdon,   
 2004,   http://numbers.computation.free.fr/Constants/Miscellaneous/zetazeros1e13-1e24.pdf.

\bibitem{Odlyzko}  
 A. M.  Odlyzko,  The $10^{21}$-st zero of the Riemann zeta function,
 www.research.att.com/$\sim$amo,  1998.  

\bibitem{Hiary}    G.  Hiary, 
Ann.  Math., 174-2 (2011)  891;   
\\
 https://people.math.osu.edu/hiary.1/outd3/out.88837796029624663862630219091085.zeros




 

\end{thebibliography}
\end{document}